\documentclass[10pt]{article}
\usepackage{verbatim}
   \newtheorem{Theorem}{Theorem}[section]

   \newtheorem{Definition}[Theorem]{Definition}
   \newtheorem{Example}[Theorem]{Example}

   \newtheorem{Conjecture}[Theorem]{Conjecture}

   \newcommand{\qed}{\hphantom{.}\hfill \rule[1pt]{8pt}{8pt}\medbreak}

   \def\ppmod#1{\allowbreak \mkern10mu ({\rm mod}\,\,#1)} 

\begin{document}
   \pagestyle{myheadings}

   \begin{center}
   \LARGE
   \textbf{Complete Mixed Doubles Round Robin Tournaments}\\
   \vspace{.5in}

   \normalsize

   David R. Berman \\
   Department of Computer Science \\
   University of North Carolina Wilmington \\
   Wilmington, NC 28403 \\
   {\texttt bermand@uncw.edu}\\

   \medskip

   Ian N. Wakeling \\
	 Qi Statistics Ltd. \\
	 Penhales House, Ruscombe \\
   Berkshire RG10 9JN, UK \\
   {\texttt ian@qistatistics.co.uk}\\
	
   \end{center}

   \vspace{.1in}


\begin{abstract}
We present a new type of tournament design that we call a complete
mixed doubles round robin tournament, CMDRR$(n,k)$, that generalizes spouse-avoiding mixed 
doubles round robin tournaments and strict Mitchell mixed doubles round robin tournaments.
We show that CMDRR$(n,k)$ exist for all allowed values of $n$ and $k$ apart from $4$ exceptions and
$31$ possible exceptions. We show that a fully resolvable CMDRR$(2n, 0)$ exists for all $n \ge 5$ 
and a fully resolvable CMDRR$(3n, n)$ exists for all $n \ge 5$ and $n$ odd. We prove a product
theorem for constructing CMDRR$(n,k)$.

\medskip
\noindent {{\it Keywords\/}: complete mixed doubles round robin tournament; 
spouse-avoiding mixed doubles round robin tournament; 
strict Mitchell mixed doubles round robin tournament; 
holey self-orthogonal latin square; fully resolvable.}

\medskip
\noindent
 {{\it AMS subject classification\/}: 05B30, 05B15.}
   \smallskip
\end{abstract}

   \section{Introduction}
   \label{s:intro}

A mixed doubles tournament is a set of games or matches between two teams, 
where each team consists of one male and one female player, as in mixed doubles tennis.  
We are concerned here with the situation in which the teams are not fixed, 
but vary throughout the tournament, unlike, say, the usual arrangement in a bridge 
tournament, where the same two players form a team in every match they play. 
Also we impose round robin properties on the tournament structure. These properties specify 
the number of times players oppose, and the number of times players of the opposite sex partner. 
The best-known type of a mixed doubles tournament in which partners are not fixed is the 
spouse-avoiding mixed doubles round robin tournament. 

A spouse-avoiding mixed doubles round robin tournament, \\ SAMDRR$(n)$, is a schedule of 
mixed doubles games for $n$ husband and wife couples. The tournament is structured so that 
spouses never play in a match together as partners or opponents. However, every man 
and woman who are not spouses are partners exactly once and opponents exactly once, 
and every pair of players of the same sex are opponents exactly once. 
Brayton, Coppersmith, and Hoffman~\cite{BCH1, BCH2} defined these tournaments in 1973 and 
showed that a SAMDRR$(n)$ exists for all $n$ except $2$, $3,$ and $6$. 
A SAMDRR$(n)$ is resolvable if the games can be arranged in rounds so that: if $n$ is even, each player plays in every round; and if $n$ is odd, each player 
except one husband and wife plays in every round. It follows that for $n$ odd, every player has exactly one bye, i.e., round they sit out. 
The total number of games is $n(n-1)/2$. 
The existence of a resolvable SAMDRR$(n)$ is equivalent to the existence of a self-orthogonal latin square of 
order $n$ with a symmetric orthogonal mate (SOLSSOM) (see \cite{IA1, FIN}). 

Recently, a new class of mixed doubles tournaments with round robin properties has been introduced and studied by 
Berman and Smith~\cite{BS1, BS2}. They are called strict Mitchell mixed doubles round robin tournaments (strict MMDRR) and 
were motivated by an article of Anderson~\cite{IA2} who describes a problem of Mitchell~\cite{M} from 
the late nineteenth century. 

\begin{Definition}
\label{strict} A strict Mitchell mixed doubles round robin tournament \\
(strict MMDRR$(n)$) is a schedule of mixed doubles games for $n$ men and $n$
women in which every man and woman partner exactly once and oppose exactly once.
Every pair of players of the same sex oppose at least once.
\end{Definition}

Note first that since every player appears in $n$ games, every player must oppose one player of the same sex exactly twice. 
Also, the number of games in a strict MMDRR$(n)$ is ${n^2}/2$. It follows that this tournament structure can be considered only 
when $n$ is even. Berman and Smith~\cite{BS1} give examples of strict MMDRR$(n)$ for 
$n = 2, 4, 6, 10, 14$, prove a product theorem, and show strict MMDRR$(n)$ exist for $n = 16k$  for $k\ge1$ and 
$n=16k+4$ for $k\ge3$. 
 
In this paper we introduce a new type of tournament called a complete mixed doubles round robin tournament that generalizes 
both SAMDRRs and strict MMDRRs. 

\begin{Definition}\label{CCMDRR} A complete mixed doubles round robin tournament \\
(CMDRR$(n, k)$) is a schedule of mixed doubles games for $n$ men and $n$
women of which $k$ men and $k$ women are spouses. Spouses never play in a match together as partners or opponents.  However, every man 
and woman who are not spouses are partners exactly once and opponents exactly once. Each player 
who has a spouse opposes every same sex player exactly once. Each player who does not have a spouse opposes 
some other same sex player who does not have a spouse exactly twice and opposes all other same sex players exactly once.
\end{Definition}

By definition every CMDRR$(n, 0)$ is a strict MMDRR$(n)$ and every CMDRR$(n, n)$ is a SAMDRR$(n)$. 
For odd $n$, CMDRR$(n, 1)$ is the closest that it is possible to come to the non-existent strict MMDRR$(n)$. The number of games in a 
CMDRR$(n, k)$ is $({n^2}-k)/2$. Players who do not have spouses are paired by repeated opposition so  $n-k$ must be even. 
We represent a CMDRR$(n, k)$ as a square matrix of order $n$ with males as row indices and females as column indices. 
The entry in position (M$i$, F$j$) is the pair (M$x$, F$y$) if and only if the game 
M$i$, F$j$ v M$x$, F$y$ is in the tournament. Each game contributes two entries, i.e., the entry in position (M$x$, F$y$) is the pair (M$i$, F$j$). 
If the CMDRR is a SAMDRR then this representation is different than the standard representation. In the standard representation, a SAMDRR$(n)$ corresponds 
to a SOLS$(n)$ with males as both row and column indices and females as entries. There is a game M$i$, F$x$ v M$j$, F$y$ if and only if the entry in position 
($i$, $j$) is $x$ and the entry in position ($j$, $i$) is $y$. This standard representation cannot be used for a CMDRR because of repeated opposition of 
same sex players. 

\begin{Example}\label{C31}  A CMDRR$(3, 1)$
\\[.05in]
\begin{tabular}{c|ccc}
    & F1   & F2   & F3  \\
\cline{1-4}
M1  &      & M2F3 & M3F2  \\
M2  & M3F3 & M3F1 & M1F2  \\
M3  & M2F2 & M1F3 & M2F1  \\
\end{tabular}
\end{Example}

From row $1$ we see that M$1$ opposes M$2$, M$3$, F$3$, F$2$, and partners F$2$, F$3$. 
From row $2$ we see that M$2$ opposes M$3$ twice and M$1$ once, opposes F$3$, F$1$, F$2$, and partners F$1$, F$2$, F$3$.
From row $3$ we see that M$3$ opposes M$2$ twice and M$1$ once, opposes F$2$, F$3$, F$1$, and partners F$1$, F$2$, F$3$. 
From the columns we see similar information about each female. The hole at position (M$1$, F$1$) indicates
that these two players are spouses. Thus, it is easy to check that the conditions for a CMDRR$(3, 1)$
are satisfied.  In future examples we will suppress the row and column headers and also the M and F in each entry. 

We next discuss resolvability for CMDRR$(n, k)$. The games must be partitioned into rounds so that each player plays in at most one game per round. 
We will call a round \textit{full} if it involves all players if $n$ is even, 
and all but $2$ players if $n$ is odd. A round that is not full is called \textit{short}. A CMDRR$(n, k)$ is called \textit{fully resolvable} if the 
games can be partitioned into rounds with at most one short round. 
The round structure is specified by a matrix of order $n$, with entries from the set \{1, \ldots, r\}, where $r$ is the number of rounds, 
and each entry appears at most one time in each row and column. The entry in cell $(i, j)$ is the round in which the game partnering M$i$ and F$j$ is played 
for non-spouses, or empty for spouses.  

\begin{Example}\label{C31b}  Resolution for the CMDRR$(3, 1)$ of Example~\ref{C31} into $4$ full rounds. 
\\[.05in]
\begin{tabular}{|ccc|}
\cline{1-3}
  & 1 & 2 \\
3 & 4 & 1 \\
4 & 2 & 3 \\
\cline{1-3}
\end{tabular}
\hspace{0.5 in}
\begin{tabular}{ccc}
round & game & byes \\
\cline{1-3}
1 & M1F2 v M2F3 & M3, F1 \\
2 & M1F3 v M3F2 & M2, F1 \\
3 & M2F1 v M3F3 & M1, F2 \\
4 & M2F2 v M3F1 & M1, F3 \\
\cline{1-3}
\end{tabular}
\end{Example}

Unfortunately, full resolvability is usually hard or impossible to come by. Alternatively we will settle for 
a partition of the games into all short rounds, all but one of equal length. Notice that every non-spouse player will have the same number of byes (say b), 
and every spouse player will have $b+1$ byes. Ideally each of the equal length short rounds should have the greatest possible number of players.

   \section{Examples}

In this section we give examples of CMDRR$(n, k)$ for $n \le 8$, and also of a CMDRR$(9, 3)$ and a CMDRR$(10, 2)$. Examples of 
SAMDRR$(n)$ can be found in~\cite{IA1}. Most of the examples were found using an Embarcadero Delphi XE program, available from the second author. 
The program fixes the partnerships and then exchanges them between games in a tabu search algorithm that seeks to optimize the 
opposition pairs incidence matrix (see~\cite{Mor}). The examples will be used in the next section as the basis of our recursive construction. 
A more extensive list of examples is available from the authors. 

The strict MMDRR$(2)$, CMDRR$(3, 1)$, and SAMDRR$(n)$ for $n=4,5,7,$ and $8$ are fully resolvable. It is easy to check by hand that the strict MMDRR$(4)$ and the CMDRR$(5,1)$ are not fully resolvable. A computer search shows that the other examples are not fully resolvable. 
A non-trivial resolution into short rounds is given when known. 
In the next section we will give general results on resolvability. 

\subsection{Tournaments with 4 players}

A SAMDRR$(2)$ does not exist. 

\begin{Example}\label{C20}  A strict MMDRR$(2)$ with repeat oppositions M$1$M$2$ and \\ F$1$F$2$. 
\\[.05in]
\begin{tabular}{|cc|}
\cline{1-2}
22 & 21  \\
12 & 11  \\
\cline{1-2}
\end{tabular}
\end{Example}
	
\subsection{Tournaments with 6 players}	

A CMDRR$(3, 1)$ was given in Example~\ref{C31}. A SAMDRR$(3)$ does not exist.
	
\subsection{Tournaments with 8 players}

A SAMDRR$(4)$ exists. A CMDRR$(4, 2)$ does not exist. 

\begin{Example}\label{C40}  A strict MMDRR$(4)$ with repeat oppositions M$1$M$2$, M$3$M$4$, \\ F$1$F$2$, and F$3$F$4$. 
\\[.05in]
\begin{tabular}{|cccc|}
\cline{1-4}
24 & 41 & 22 & 33 \\
32 & 13 & 44 & 11 \\
43 & 21 & 14 & 42 \\
12 & 34 & 31 & 23 \\
\cline{1-4}
\end{tabular}
\end{Example}

\subsection{Tournaments with 10 players}

A SAMDRR$(5)$ exists.

\begin{Example}\label{C51}  A CMDRR$(5, 1)$ with spouse pair M$1$F$1$ and repeat oppositions 
M$2$M$3$, M$4$M$5$, F$2$F$3$, and F$4$F$5$. 
\\[.05in]
\begin{tabular}{|ccccc|}
\cline{1-5}
   & 55 & 22 & 43 & 34 \\
54 & 13 & 32 & 35 & 41 \\
42 & 23 & 51 & 15 & 24 \\
25 & 31 & 14 & 52 & 53 \\
33 & 44 & 45 & 21 & 12 \\
\cline{1-5}
\end{tabular}
\end{Example}

\begin{Conjecture} A CMDRR$(5, 3)$ does not exist.
\end{Conjecture}

\subsection{Tournaments with 12 players}

A SAMDRR$(6)$ does not exist.

\newpage
\begin{Example}\label{C60}  A strict MMDRR$(6)$ with repeat oppositions M$1$M$2$, M$3$M$4$, \\ M$5$M$6$, 
F$1$F$2$, F$3$F$4$, and F$5$F$6$. 
\\[.05in]
\begin{tabular}{|cccccc|}
\cline{1-6}
55 & 63 & 24 & 42 & 26 & 31 \\
44 & 36 & 61 & 13 & 52 & 15 \\
16 & 41 & 45 & 53 & 64 & 22 \\
32 & 14 & 56 & 21 & 33 & 65 \\
62 & 25 & 34 & 66 & 11 & 43 \\
23 & 51 & 12 & 35 & 46 & 54 \\
\cline{1-6}
\end{tabular}
\\[.05in]
Resolution with short rounds $1$-$9$. Each round has $2$ games and each player has $3$ byes.
\\[.05in]
\begin{tabular}{|cccccc|}
\cline{1-6}
8 & 1 & 6 & 7 & 5 & 3 \\
1 & 8 & 7 & 6 & 3 & 5 \\
3 & 6 & 9 & 2 & 4 & 8 \\
6 & 7 & 4 & 1 & 9 & 2 \\
5 & 3 & 2 & 9 & 8 & 4 \\
7 & 5 & 1 & 4 & 2 & 9 \\
\cline{1-6}
\end{tabular}
\end{Example}

\begin{Conjecture} A CMDRR$(6, 2)$ does not exist.
\end{Conjecture}

\begin{Example}\label{C64}  A CMDRR$(6, 4)$ with spouse pairs M$1$F$1$, M$2$F$2$, M$3$F$3$, and M$4$F$4$, 
and repeat oppositions M$5$M$6$ and F$5$F$6$.
\\[.05in]
\begin{tabular}{|cccccc|}
\cline{1-6}
   & 34 & 52 & 66 & 43 & 25 \\
54 &    & 41 & 35 & 16 & 63 \\
46 & 61 &    & 12 & 24 & 55 \\
23 & 56 & 15 &    & 62 & 31 \\
65 & 13 & 64 & 21 & 36 & 42 \\
32 & 45 & 26 & 53 & 51 & 14 \\
\cline{1-6}
\end{tabular}
\end{Example}

\subsection{Tournaments with 14 players}

A SAMDRR$(7)$ exists.

\begin{Example}\label{C71}   A CMDRR$(7, 1)$ with spouse pair M$1$F$1$ and repeat oppositions 
M$2$M$3$, M$4$M$5$, M$6$M$7$, F$2$F$3$, F$4$F$5$, and F$6$F$7$. 
\\[.05in]
\begin{tabular}{|ccccccc|}
\cline{1-7}
   & 75 & 32 & 43 & 64 & 27 & 56 \\
53 & 31 & 65 & 37 & 74 & 42 & 16 \\
22 & 13 & 77 & 66 & 41 & 55 & 24 \\
35 & 26 & 14 & 51 & 57 & 73 & 62 \\
44 & 63 & 21 & 72 & 36 & 17 & 45 \\
76 & 47 & 52 & 15 & 23 & 34 & 71 \\
67 & 54 & 46 & 25 & 12 & 61 & 33 \\
\cline{1-7}
\end{tabular}
\\[.05in]
\newpage\noindent
Resolution with short rounds $1-12$. Each round has $2$ games, each player 
except M$1$ and F$1$ has $5$ byes, and each of M$1$ and F$1$ have $6$ byes.
\\[.05in]
\begin{tabular}{|ccccccc|}
\cline{1-7}
  & 8 & 9 & 2 & 5 & 6 & 3 \\
7 & 2 & 3 & 12 & 1 & 10 & 4 \\
2 & 9 & 10 & 8 & 11 & 4 & 12 \\
11 & 10 & 2 & 6 & 9 & 5 & 7 \\
6 & 1 & 7 & 11 & 4 & 3 & 9 \\
12 & 7 & 1 & 5 & 3 & 8 & 4 \\
4 & 11 & 5 & 1 & 8 & 12 & 10 \\
\cline{1-7}
\end{tabular}
\end{Example}

\begin{Example}\label{C73}   A CMDRR$(7, 3)$ with spouse pairs M$1$F$1$, M$2$F$2$, and M$3$F$3$, 
and repeat oppositions M$4$M$5$, M$6$M$7$, F$4$F$5$, and F$6$F$7$. 
\\[.05in]
\begin{tabular}{|ccccccc|}
\cline{1-7}
   & 55 & 62 & 26 & 43 & 37 & 74 \\
57 &    & 41 & 73 & 36 & 14 & 65 \\
64 & 47 &    & 52 & 71 & 25 & 16 \\
23 & 76 & 15 & 67 & 54 & 51 & 32 \\
46 & 34 & 77 & 45 & 12 & 63 & 21 \\
72 & 13 & 56 & 31 & 27 & 75 & 44 \\
35 & 61 & 24 & 17 & 66 & 42 & 53 \\
\cline{1-7}
\end{tabular}
\end{Example}

\begin{Example}\label{C75}   A CMDRR$(7, 5)$ with spouse pairs M$1$F$1$, M$2$F$2$, M$3$F$3$,
M$4$F$4$, and M$5$F$5$, and repeat oppositions M$6$M$7$ and F$6$F$7$. 
\\[.05in]
\begin{tabular}{|ccccccc|}
\cline{1-7}
   & 54 & 42 & 35 & 66 & 27 & 73 \\
45 &    & 61 & 53 & 37 & 74 & 16 \\
52 & 76 &    & 67 & 14 & 41 & 25 \\
36 & 13 & 75 &    & 21 & 57 & 62 \\
77 & 31 & 24 & 12 &    & 63 & 46 \\
23 & 47 & 56 & 71 & 72 & 15 & 34 \\
64 & 65 & 17 & 26 & 43 & 32 & 51 \\
\cline{1-7}
\end{tabular}
\end{Example}

\subsection{Tournaments with 16 players}

A SAMDRR$(8)$ exists.

\begin{Example}\label{C80}   A strict MMDRR$(8)$ with repeat oppositions 
M$1$M$2$, \\ M$3$M$4$, M$5$M$6$, M$7$M$8$, F$1$F$2$, F$3$F$4$, F$5$F$6$, and F$7$F$8$. 
\\[.05in]
\begin{tabular}{|cccccccc|}
\cline{1-8}
84 & 45 & 52 & 26 & 67 & 73 & 38 & 21 \\
18 & 37 & 65 & 42 & 76 & 14 & 51 & 83 \\
63 & 71 & 44 & 85 & 56 & 48 & 22 & 17 \\
55 & 24 & 87 & 33 & 12 & 61 & 78 & 36 \\
27 & 13 & 74 & 68 & 41 & 35 & 86 & 62 \\
46 & 58 & 31 & 77 & 23 & 82 & 15 & 54 \\
32 & 81 & 16 & 53 & 88 & 25 & 64 & 47 \\
72 & 66 & 28 & 11 & 34 & 57 & 43 & 75\\
\cline{1-8}
\end{tabular}
\\[.05in]
Resolution with short rounds $1-10$, each with $3$ games,  and short round $11$ with $2$ games. 
Each player has $3$ byes. 
\\[.05in]
\begin{tabular}{|cccccccc|}
\cline{1-8}
7 & 9 & 8 & 11 & 3 & 10 & 4 & 2 \\
2 & 10 & 1 & 6 & 7 & 11 & 5 & 9 \\
6 & 1 & 5 & 8 & 2 & 3 & 10 & 4 \\
11 & 6 & 2 & 5 & 9 & 4 & 8 & 3 \\
5 & 8 & 4 & 10 & 11 & 2 & 1 & 7 \\
4 & 7 & 6 & 1 & 9 & 5 & 3 & 10 \\
1 & 3 & 10 & 4 & 6 & 7 & 9 & 8 \\
3 & 5 & 9 & 7 & 8 & 1 & 2 & 6\\
\cline{1-8}
\end{tabular}
\end{Example}

\begin{Example}\label{C82}  A CMDRR$(8, 2)$ with spouse pairs M$1$F$1$ and M$2$F$2$, and repeat oppositions 
M$3$M$8$, M$4$M$5$, M$6$M$7$, F$3$F$7$, F$4$F$8$, and F$5$F$6$. 
\\[.05in]
\begin{tabular}{|cccccccc|}
\cline{1-8}
   & 65 & 87 & 38 & 74 & 23 & 56 & 42 \\
75 &    & 16 & 47 & 58 & 64 & 81 & 33 \\
46 & 83 & 28 & 51 & 67 & 85 & 72 & 14 \\
62 & 18 & 55 & 73 & 86 & 31 & 24 & 57 \\
34 & 76 & 61 & 82 & 43 & 17 & 48 & 25 \\
53 & 41 & 77 & 26 & 12 & 78 & 35 & 84 \\
88 & 37 & 44 & 15 & 21 & 52 & 63 & 66 \\
27 & 54 & 32 & 68 & 36 & 45 & 13 & 71\\
\cline{1-8}
\end{tabular}

\end{Example}

\begin{Example}\label{C84}  A CMDRR$(8, 4)$ with spouse pairs M$1$F$1$, M$2$F$2$, 
M$3$F$3$ and M$4$F$4$ and repeat oppositions 
M$5$M$7$, M$6$M$8$, F$5$F$7$, and F$6$F$8$. 
\\[.05in]
\begin{tabular}{|cccccccc|}
\cline{1-8}
   & 58 & 42 & 83 & 66 & 34 & 25 & 77 \\
73 &    & 55 & 61 & 17 & 48 & 36 & 84 \\
88 & 64 &    & 16 & 72 & 27 & 51 & 45 \\
52 & 13 & 67 &    & 38 & 71 & 85 & 26 \\
37 & 41 & 86 & 75 & 23 & 68 & 74 & 12 \\
24 & 87 & 78 & 32 & 81 & 15 & 43 & 56 \\
46 & 35 & 21 & 57 & 54 & 82 & 18 & 63 \\
65 & 76 & 14 & 28 & 47 & 53 & 62 & 31\\
\cline{1-8}
\end{tabular}
\end{Example}

\begin{Example}\label{C86}  A CMDRR$(8, 6)$ with spouse pairs M$1$F$1$, M$2$F$2$, 
M$3$F$3$, M$4$F$4$, M$5$F$5$, and M$6$F$6$, and repeat oppositions 
M$7$M$8$  and F$7$F$8$. 
\\[.05in]
\begin{tabular}{|cccccccc|}
\cline{1-8}
   & 64 & 56 & 38 & 47 & 72 & 83 & 25 \\
43 &    & 74 & 86 & 18 & 35 & 61 & 57 \\
68 & 75 &    & 51 & 26 & 87 & 42 & 14 \\
76 & 37 & 21 &    & 63 & 58 & 15 & 82 \\
34 & 81 & 62 & 77 &    & 13 & 28 & 46 \\
27 & 53 & 45 & 12 & 84 &    & 78 & 31 \\
85 & 16 & 88 & 23 & 32 & 41 & 54 & 67 \\
52 & 48 & 17 & 65 & 71 & 24 & 36 & 73\\
\cline{1-8}
\end{tabular}
\end{Example}

\subsection{Tournaments with more than 16 players}

\begin{Example}\label{C93}  A CMDRR$(9, 3)$ with spouse pairs M$1$F$1$, M$2$F$2$, and\\
M$3$F$3$, and repeat oppositions M$4$M$5$, F$4$F$5$, M$6$M$7$, F$6$F$7$, M$8$M$9$, and F$8$F$9$. 
\\[.05in]
\begin{tabular}{|ccccccccc|}
\cline{1-9}
   & 47 & 29 & 66 & 53 & 82 & 38 & 75 & 94 \\
77 &    & 68 & 35 & 46 & 91 & 54 & 89 & 13 \\
95 & 61 &    & 52 & 24 & 78 & 49 & 17 & 86 \\
63 & 58 & 56 & 81 & 74 & 25 & 12 & 99 & 37 \\
88 & 34 & 15 & 27 & 69 & 43 & 96 & 42 & 71 \\
32 & 79 & 41 & 98 & 87 & 14 & 76 & 23 & 55 \\
59 & 93 & 84 & 45 & 18 & 67 & 21 & 36 & 62 \\
44 & 16 & 97 & 73 & 92 & 39 & 65 & 51 & 28 \\
26 & 85 & 72 & 19 & 31 & 57 & 83 & 64 & 48 \\
\cline{1-9}
\end{tabular}
\\[.05in]
Resolution with short rounds $1-13$, each with $3$ games. 
\\[.05in]
\begin{tabular}{|ccccccccc|}
\cline{1-9}
   & 3  & 4  & 9  & 2  & 12 & 1  & 11 & 13 \\
6  &    & 3  & 12 & 5  & 2  & 7  & 10 & 4  \\
10 & 5  &    & 11 & 12 & 8  & 2  & 1  & 6  \\
7  & 13 & 1  & 8  & 4  & 5  & 3  & 12 & 2  \\
9  & 11 & 2  & 7  & 8  & 1  & 4  & 13 & 3  \\
5  & 1  & 7  & 6  & 13 & 9  & 10 & 3  & 8  \\
3  & 9  & 5  & 4  & 11 & 10 & 6  & 8  & 1  \\
8  & 12 & 11 & 5  & 7  & 6  & 13 & 9  & 10 \\
2  & 7  & 9  & 13 & 10 & 4  & 11 & 6  & 12 \\
\cline{1-9}
\end{tabular}
\end{Example}

\begin{Example}\label{C102}  A CMDRR$(10, 2)$ with spouse pairs M$1$F$1$ and M$2$F$2$, 
and repeat oppositions M$3$M$4$, F$3$F$4$, M$5$M$6$, F$5$F$6$, M$7$M$8$, F$7$F$8$, M$9$M$0$, and F$9$F$0$. 
\\[.05in]
\begin{tabular}{|cccccccccc|}
\cline{1-10}
   & 40 & 74 & 92 & 66 & 85 & 29 & 37 & 53 & 08 \\
48 &    & 56 & 39 & 00 & 71 & 84 & 95 & 17 & 63 \\
90 & 05 & 82 & 51 & 73 & 67 & 18 & 49 & 24 & 46 \\
69 & 86 & 04 & 75 & 57 & 30 & 93 & 21 & 38 & 12 \\
34 & 68 & 19 & 60 & 81 & 23 & 45 & 76 & 02 & 97 \\
03 & 77 & 20 & 88 & 99 & 15 & 36 & 52 & 41 & 54 \\
26 & 91 & 35 & 13 & 44 & 58 & 62 & 07 & 80 & 89 \\
55 & 33 & 98 & 27 & 16 & 42 & 01 & 64 & 70 & 79 \\
72 & 14 & 47 & 06 & 28 & 09 & 50 & 83 & 65 & 31 \\
87 & 59 & 61 & 43 & 32 & 94 & 78 & 10 & 96 & 25 \\
\cline{1-10}
\end{tabular}
\\[.05in]
\newpage\noindent
Resolution with short rounds $1-16$, each with $3$ games, and short round $17$ with $1$ game. 
\\[.05in]
\begin{tabular}{|cccccccccc|}
\cline{1-10}
   & 7  & 2  & 5  & 15 & 3  & 9  & 6  & 1  & 16 \\
14 &    & 16 & 3  & 6  & 12 & 13 & 2  & 9  & 8  \\
4  & 11 & 12 & 15 & 7  & 14 & 6  & 13 & 3  & 1  \\
11 & 6  & 4  & 16 & 10 & 1  & 3  & 14 & 13 & 7  \\
15 & 4  & 1  & 9  & 8  & 16 & 10 & 5  & 2  & 11 \\
13 & 1  & 8  & 10 & 12 & 15 & 14 & 4  & 11 & 9  \\
12 & 10 & 7  & 2  & 16 & 5  & 1  & 15 & 17 & 14 \\
8  & 12 & 9  & 13 & 3  & 6  & 5  & 10 & 14 & 17 \\
10 & 5  & 3  & 8  & 2  & 7  & 11 & 9  & 12 & 4  \\
5  & 2  & 13 & 4  & 11 & 8  & 15 & 16 & 7  & 6  \\
\cline{1-10}
\end{tabular}
\end{Example}

   \section{Recursive Construction}
	
In this section we present a recursive construction using holey SOLS and use it to show the existence of 	
CMDRR$(n,k)$ for all allowed values of $n$ and $k$, apart from $4$ exceptions and $31$ possible exceptions.  
We show that a fully resolvable CMDRR$(2n, 0)$ exists for all $n \ge 5$ 
and a fully resolvable CMDRR$(3n, n)$ exists for all $n \ge 5$ and $n$ odd.

For completeness we include the definition of a holey SOLS (see~\cite{FIN}).

\begin{Definition}\label{HSOLS} 
A holey SOLS (or frame SOLS) is a self-orthogonal latin square of order $n$ with $n_i$ missing sub-SOLS (holes)
of order $h_i$ ($1 \le i \le k$), which are disjoint and spanning (that is $\sum_{1 \le i \le k}n_ih_i = n$). 
It is denoted by HSOLS$(h_1^{n_1} \ldots h_k^{n_k})$ where $h_1^{n_1} \ldots h_k^{n_k}$ is the type of the HSOLS. 
\end{Definition}
 
Suppose an HSOLS exists and CMDRR exist for each hole size. Then we can fill in the holes with the CMDRRs to get a 
new CMDRR. The details follow. For convenience we will assume that a CMDRR$(1, 1)$ exists with spouse pair M$1$F$1$ and no games. 

\begin{Theorem}\label{cons}
Suppose an HSOLS$(n)$ of type $h_1^{n_1} \ldots h_k^{n_k}$ exists and for each $h_i$ there exist 
CMDRR$(h_i, m_{i1}), \ldots, $ CMDRR$(h_i, m_{i{n_i}})$. Then there exists a CMDRR$(n, s)$ where 
$s = \sum_{1 \le i \le k}\sum_{1 \le j \le n_i}m_{ij}$.
\end{Theorem}
\begin{Proof}
By possibly relabeling players we can assume that the HSOLS is block diagonal. By construction, spouse pairs will always have the form M$i$F$i$. 

Use the standard SOLS representation for a SAMDRR to identify games for all entries of the HSOLS that are not in a hole. 
Thus every entry $(i,j)$ not in a hole will 
contribute the game M$i$F$(i,j)$ v M$j$F$(j,i)$. By definition of an HSOLS these games satisfy the conditions that every pair of opposite sex players 
have partnered and opposed at most once, every pair of same sex players have opposed at most once, and no spouses have played in a game together. 

Each missing sub-SOLS (hole) corresponds to a set S of consecutive integers which are the indices and also the missing entries. Use a translation of an appropriate sized
CMDRR and identify games using the representation introduced for CMDRRs. Assume the translation is by $t$. 
Then a non-empty entry $(i,j)$ of the CMDRR will contribute the game M$i'$M$j'$ v F$x'$F$y'$ where 
the $(i,j)$ entry of the CMDRR is $(x,y)$ and each primed symbol is the corresponding unprimed symbol plus $t$. By definition of an HSOLS, the players in a hole will 
not be involved in any other common games. 

Taking all the games identified by the two different processes described above will produce the required CMDRR. 
As every pair of players is either in a hole or not in any hole, the conditions for a CMDRR are met.  \qed
\end{Proof}

\begin{Example}\label{C11} A CMDRR$(11, 7)$  can be constructed from the \\ HSOLS$(1^63^12^1)$ given below (see~\cite{FIN, Z}). Use the 
CMDRR$(3, 1)$ of Example~\ref{C31} to fill the hole of size three and the  strict MMDRR$(2)$ of Example~\ref{C20} to fill the hole of size two. 
The spouse pairs are M$1$F$1, \ldots, $M$7$F$7$. The repeat oppositions are M$8$M$9$, F$8$F$9$, M$10$M$11$, and F$10$F$11$.
\\[.05in]
\begin{tabular}{|ccccccccccc|}
\cline{1-11}
. & 7 & 5 & 2 & 10 & 9 & 4 & 11 & 3 & 6 & 8 \\
6 & . & 8 & 7 & 3 & 10 & 1 & 5 & 11 & 4 & 9 \\
9 & 4 & . & 10 & 8 & 1 & 11 & 2 & 6 & 5 & 7 \\
5 & 11 & 9 & . & 7 & 2 & 6 & 10 & 1 & 8 & 3 \\
7 & 6 & 11 & 3 & . & 8 & 2 & 4 & 10 & 9 & 1 \\
11 & 8 & 4 & 9 & 1 & . & 10 & 3 & 5 & 7 & 2 \\
3 & 10 & 6 & 1 & 11 & 5 & . & . & . & 2 & 4 \\
4 & 1 & 10 & 6 & 2 & 11 & . & . & . & 3 & 5 \\
10 & 5 & 2 & 11 & 4 & 3 & . & . & . & 1 & 6 \\
8 & 9 & 7 & 5 & 6 & 4 & 3 & 1 & 2 & . & . \\
2 & 3 & 1 & 8 & 9 & 7 & 5 & 6 & 4 & . & . \\
\cline{1-11}
\end{tabular}
\end{Example}

In order to create CMDRR using Theorem~\ref{cons} we need a supply of HSOLS. The following theorems (see~\cite{FIN, HZ, XZZ}) give just the 
supply we need.

\begin{Theorem}\label{HSOLSanb1}
For $n \ge 4$ and $a \ge 2$, an HSOLS$(a^nb^1)$ exists if $0 \le b \le a(n-1)/2$ with 
possible exceptions for $n \in \{6, 14, 18, 22\}$ and $b = a(n-1)/2$.
\end{Theorem}

\begin{Definition}
An incomplete SOLS is a self-orthogonal latin square of order $n$ missing a sub-SOLS of 
order $k$, denoted by ISOLS$(n,k)$. An ISOLS$(n,k)$ is equivalent to an HSOLS$(1^{n-k}k^1)$. (see~\cite{FIN})
\end{Definition}

\begin{Theorem}\label{ISOLS}
There exists an ISOLS$(n,k)$ for all values of $n$ and $k$ satisfying $n \ge 3k+1$, 
except for $(n,k) = (6,1), (8,2)$ and possibly excepting $n = 3k+2$ and $k \in \{6,8,10\}$.
\end{Theorem}

We can now show that CMDRR$(n,k)$ exist for all but a finite number of possible exceptions. 

\begin{Theorem}\label{ge32}
There exists a CMDRR$(n,k)$ for each $n \ge 32, k \le n, $ and $n-k$ even. 
\end{Theorem}
\begin{Proof}
By Theorem~\ref{HSOLSanb1} an HSOLS$(a^nb^1)$ exists for $a=8, n \ge 4, $ and $0 \le b \le 7$. 
Each hole of size $8$ can be filled with any one of the CMDRR$(8,0)$ (Example~\ref{C80}), CMDRR$(8,2)$ (Example~\ref{C82}), 
CMDRR$(8,4)$ \\ (Example~\ref{C84}), CMDRR$(8,6)$ (Example~\ref{C86}), or a SAMDRR$(8)$. For $2 \le b \le 7$, 
the hole of size $b$ can be filled with one of the CMDRR$(b,i)$ given in Examples~\ref{C20}--~\ref{C75} or 
the SAMDRR$(b)$ for $b \ne 2,3,6$. By direct observation and using Theorem~\ref{HSOLSanb1} we have the stated
CMDRRs, except for CMDRR$(n,n)$ when $n \ge 32$ and $n \equiv 2, 3,$ or $6 \ppmod 8$. But we know that SAMDRR$(n)$ exist for these sizes. \qed
\end{Proof}

Finally we can use Theorems~\ref{HSOLSanb1} and~\ref{ISOLS} to handle special cases for $n < 32$ and give our main result. 

\begin{Theorem}
There exists a CMDRR$(n,k)$ for each $n \ge 2, k \le n, $ and \\ $n-k$ even, except for $(n,k) = (2,2),$ $(3,3),$ $(4,2),$ $(6,6)$ 
and possibly excepting the following $31$ values: $(n,k) = (5,3),$ $(6,2),$ 
 $(12,2),$ $(12,6),$ $(12,8),$ $(13,3),$ $(13,7),$ $\!(14,2),$ $(14,6),$ $(15,3),$ $(15,7),$ $(15,9),$ $(16,2),$
 $\!(16, 10),$ $\!(17,3),$ $\!(17,7),$ $\!(17,11),$ $\!(18,2),$ $(18,10),$ $(19,3),$ $(19,11),$ $(20,2),$
 $(20,14),$ $(21,3),$ $(21, 11),$ $(22,2),$ $(22,14),$ $(23,15),$ $(24,2),$ $(24,14),$
 $(25,15).$ 
\end{Theorem}
\begin{Proof}
By Theorem~\ref{ge32}, a CMDRR$(n,k)$ exists for $n \ge 32$. The table below shows a construction that can be used for each value
of $n < 32$ and compatible $k$, except for the values listed above.  
\\[.05in]
\begin{tabular}{lll}
n & k & construction  \\ \cline{1-3}
  &   &  \\

9 & 1 & HSOLS$(2^41^1)$ \\
9 & 3 & Example \ref{C93} \\
9 & 5 & HSOLS$(2^21^5)$ Lemma 2.2 of~\cite{XZZ} \\
9 & 7 & ISOLS$(9,2)$ \\

10 & 0 & HSOLS$(2^5)$ \\
10 & 2 & Example \ref{C102} \\
10 & 4 & HSOLS$(2^31^4)$ Lemma 2.2 of~\cite{XZZ} \\
10 & 6 & HSOLS$(2^21^6)$ Lemma 2.2 of~\cite{XZZ} \\
10 & 8 & ISOLS$(10,2)$ \\

11 & 1 & HSOLS$(2^51^1)$ \\
11 & 3 & HSOLS$(2^41^3)$ Lemma 2.2 of~\cite{XZZ} \\
11 & 5 & HSOLS$(2^31^5)$ Lemma 2.2 of~\cite{XZZ} \\
11 & 7 & HSOLS$(1^63^12^1)$ Example \ref{C11} \\
11 & 9 & ISOLS$(11,2)$ \\

12 & 0 & HSOLS$(2^6)$ \\
12 & 4 & HSOLS$(3^4)$ \\
12 & 10 & ISOLS$(12,2)$ \\

13 & 1 & HSOLS$(2^53^1)$ \\
13 & 5 & HSOLS$(3^41^1)$ \\
13 & 9 & ISOLS$(13,4)$ \\
13 & 11 & ISOLS$(13,2)$ \\

14 & 0 & HSOLS$(2^7)$ \\
14 & 4 & HSOLS$(3^42^1)$ \\
14 & 8 & HSOLS$(1^84^12^1)$~\cite{ABZZ} \\
14 & 10,14 & ISOLS$(14,4)$ \\
14 & 12 & ISOLS$(14,2)$ \\

15 & 1 & HSOLS$(2^63^1)$ \\
15 & 5 & HSOLS$(3^5)$ \\
15 & 11 & ISOLS$(15,4)$ \\
15 & 13 & ISOLS$(15,2)$ \\

16 & 0,4,8,12,16 & HSOLS$(4^4)$ \\
16 & 6 & HSOLS$(3^51^1)$ \\
16 & 14 & ISOLS$(16,2)$ \\

17 & 1,5,9,13,17 & HSOLS$(4^41^1)$ \\
17 & 15 & ISOLS$(17,2)$ \\ 

18 & 0,4,8,12,16 & HSOLS$(4^42^1)$ \\
18 & 6 & HSOLS$(3^6)$ \\
18 & 14,18 & ISOLS$(18,4)$ \\ 

19 & 1,5,9,13,17 & HSOLS$(4^43^1)$ \\
19 & 7 & HSOLS$(3^61^1)$ \\
19 & 15,19 & ISOLS$(19,4)$ \\ 

20 & 0,4,8,12,16,20 & HSOLS$(4^5)$ \\
20 & 6,10 & HSOLS$(3^55^1)$ \\
20 & 18 & ISOLS$(20,2)$ \\ 

\end{tabular}
\newpage
\begin{tabular}{lll}

21 & 1,5,9,13,17,21 & HSOLS$(4^51^1)$ \\
21 & 7 & HSOLS$(3^7)$ \\
21 & 15,19 & ISOLS$(21,6)$ \\

22 & 0,4,8,12,16,20 & HSOLS$(4^52^1)$ \\
22 & 6,10 & HSOLS$(3^64^1)$ \\
22 & 16,18,20,22 & ISOLS$(22,7)$ \\ 

23 & 1,5,9,13,17,21 & HSOLS$(4^53^1)$ \\
23 & 3 & HSOLS$(2^87^1)$ \\
23 & 7,11 & HSOLS$(3^65^1)$ \\ 
23 & 17,19,21,23 & ISOLS$(23,7)$ \\ 

24 & 0,4,8,12,16 & HSOLS$(6^4)$ \\
24 & 6,10 & HSOLS$(3^66^1)$ \\
24 & 18,20,22,24 & ISOLS(24,7) \\ 

25 & 1,3,5,7,9,11,13 & HSOLS$(6^41^1)$ \\
25 & 17,19,21,23,25 & ISOLS(25,8) \\ 

26 & 0,4,8,12,16,20,24 & HSOLS$(4^56^1)$ \\
26 & 2 & HSOLS$(2^98^1)$ \\
26 & 6,10,14,18,22,26 & HSOLS$(5^51^1)$ \\ 

27 & 1,3,5,7,9,11,13,15,17,19,21,23,25,27 & HSOLS$(4^57^1)$ \\ 

28 & 0,2 & HSOLS$(2^{10}8^1)$ \\
28 & 4,6,8,10,12,14,16,18,20,22,24,26,28 & HSOLS$(7^4)$ \\ 

29 & 1,3 & HSOLS$(2^{11}7^1)$ \\
29 & 5,7,9,11,13,15,17,19,21,23,25,27,29 & HSOLS$(7^41^1)$ \\ 

30 & 0,2 & HSOLS$(2^{11}8^1)$ \\
30 & 4,6,8,10,12,14,16,18,20,22,24,26,28 & HSOLS$(7^42^1)$ \\ 

31 & 1,3 & HSOLS$(2^{12}7^1)$ \\
31 & 5,7,9,11,13,15,17,19,21,23,25,27,29 & HSOLS$(7^43^1)$ \\
\end{tabular}
\\[.05in]
\qed
\end{Proof}

Resolvability of a CMDRR is more difficult to ensure. 
By filling holes in an HSOLSSOM we can construct resolvable CMDRR. For completeness we include the definition of a holey SOLSSOM (see~\cite{FIN}).

\begin{Definition}\label{HSOLSSOM} 
A holey SOLSSOM (or frame SOLSSOM) is a holey self-orthogonal latin square S of order $n$ and type $h_1^{n_1} \ldots h_k^{n_k}$,
together with a symmetric partitioned latin square M of order $n$ and type $h_1^{n_1} \ldots h_k^{n_k}$, satisfying the property
that when superimposed, the ordered pairs are exactly those pairs of symbols that are from different holes. A holey SOLSSOM with this
structure is denoted by  HSOLSSOM$(h_1^{n_1} \ldots h_k^{n_k})$, where $h_1^{n_1} \ldots h_k^{n_k}$ is the type of the HSOLSSOM.
\end{Definition}

Using the next Theorem~\cite{FIN, GA} we can construct resolvable CMDRR. 

\begin{Theorem}
An HSOLSSOM$(2^n)$ exists for all $n \ge 5$ and an \\ HSOLSSOM$(3^n)$ exists for all odd values of $n$ with $n \ge 5$.
\end{Theorem}

\begin{Theorem}
A fully resolvable strict MMDRR$(2n)$ exists for all $n \ge 5$.
\end{Theorem}
\begin{Proof}
Begin with an HSOLSSOM$(2^n)$ and convert this to a resolvable mixed doubles round robin with $2n$ rounds of play. 
Pairs of rounds will be missing all four of the players from one of the $n$ holes. Simply fill these holes with
a strict MMDRR$(2)$ constructed on the corresponding four players, thus completing the $2n$ rounds. \qed
\end{Proof}

\begin{Example}\label{x} Lemma 2.1.2 of Bennett and Zhu~\cite{BZ} gives both an example of a holey Steiner pentagon system (HSPS) of type $2^6$
and also its equivalent HSOLSSOM$(2^6)$. We rearrange the latter to make the holes block diagonal. 
\\[.05in]
\begin{tabular}{|cccccccccccc|}
\cline{1-12}
.  & .  & 6  & 10 & 12 & 11 & 5  & 3  & 4  & 8  & 9  &  7 \\
.  & .  & 10 & 11 & 8  & 9  & 4  & 6  & 5  & 12 & 7  &  3 \\
12 & 8  & .  & .  & 9  & 7  & 10 & 1  & 11 & 6  & 5  &  2 \\
6  & 5  & .  & .  & 7  & 12 & 2  & 9  & 1  & 11 & 8  &  10 \\
10 & 4  & 12 & 9  & .  & .  & 1  & 11 & 2  & 7  & 3  &  8 \\
4  & 11 & 1  & 8  & .  & .  & 12 & 2  & 7  & 3  & 10 &  9 \\
9  & 10 & 11 & 12 & 4  & 3  & .  & .  & 6  & 5  & 2  &  1 \\
11 & 3  & 9  & 6  & 2  & 10 & .  & .  & 12 & 1  & 4  &  5 \\
7  & 12 & 8  & 5  & 3  & 2  & 11 & 4  & .  & .  & 1  &  6 \\
5  & 7  & 2  & 1  & 11 & 8  & 3  & 12 & .  & .  & 6  &  4 \\
8  & 6  & 7  & 2  & 10 & 1  & 9  & 5  & 3  & 4  & .  &  . \\
3  & 9  & 5  & 7  & 1  & 4  & 6  & 10 & 8  & 2  & .  &  . \\
\cline{1-12}
\end{tabular}
\\[.05in]
\begin{tabular}{|cccccccccccc|}
\cline{1-12}
.  & .  & 7  & 11 & 8  & 10 & 4  & 9  & 5  & 12 & 3  &  6 \\
.  & .  & 6  & 8  & 11 & 7  & 12 & 10 & 3  & 4  & 9  &  5 \\
7  & 6  & .  & .  & 2  & 12 & 5  & 11 & 1  & 8  & 10 &  9 \\
11 & 8  & .  & .  & 1  & 9  & 10 & 12 & 7  & 6  & 5  &  2 \\
8  & 11 & 2  & 1  & .  & .  & 9  & 4  & 12 & 3  & 7  &  10 \\
10 & 7  & 12 & 9  & .  & .  & 1  & 3  & 11 & 2  & 4  &  8 \\
4  & 12 & 5  & 10 & 9  & 1  & .  & .  & 2  & 11 & 6  &  3 \\
9  & 10 & 11 & 12 & 4  & 3  & .  & .  & 6  & 5  & 2  &  1 \\
5  & 3  & 1  & 7  & 12 & 11 & 2  & 6  & .  & .  & 8  &  4 \\
12 & 4  & 8  & 6  & 3  & 2  & 11 & 5  & .  & .  & 1  &  7 \\
3  & 9  & 10 & 5  & 7  & 4  & 6  & 2  & 8  & 1  & .  &  . \\
6  & 5  & 9  & 2  & 10 & 8  & 3  & 1  & 4  & 7  & .  &  . \\
\cline{1-12}
\end{tabular}
\\[.05in]
\newpage\noindent
The HSOLSSOM is converted to a mixed doubles tournament and filled to produce a fully resolvable strict MMDRR$(12)$ with $12$ rounds. 
\\[.05in]
{\footnotesize
\begin{tabular}{cccccccc}
R1	& M01 F01 v M02 F02	& M03 F11 v M09 F08	& M04 F07 v M05 F09	&\\& M06 F12 v M07 F03	& M08 F05 v M12 F10	& M10 F06 v M11 F04 \\
R2	& M01 F02 v M02 F01	& M03 F09 v M05 F12	& M04 F10 v M12 F07	&\\& M06 F03 v M10 F08	& M07 F06 v M09 F11	& M08 F04 v M11 F05 \\
R3	& M01 F09 v M11 F08	& M02 F05 v M09 F12	& M03 F03 v M04 F04	&\\& M05 F07 v M10 F11	& M06 F02 v M08 F10	& M07 F01 v M12 F06 \\
R4	& M01 F05 v M07 F09	& M02 F12 v M10 F07	& M03 F04 v M04 F03	&\\& M05 F11 v M08 F02	& M06 F10 v M11 F01	& M09 F06 v M12 F08 \\
R5	& M01 F04 v M09 F07	& M02 F03 v M12 F09	& M03 F10 v M07 F11	&\\& M04 F08 v M11 F02	& M05 F05 v M06 F06	& M08 F01 v M10 F12 \\
R6	& M01 F07 v M12 F03	& M02 F10 v M03 F08	& M04 F11 v M10 F01	&\\& M05 F06 v M06 F05	& M07 F02 v M11 F09	& M08 F12 v M09 F04 \\
R7	& M01 F06 v M03 F12	& M02 F09 v M06 F11	& M04 F01 v M09 F05	&\\& M05 F03 v M11 F10	& M07 F07 v M08 F08	& M10 F04 v M12 F02 \\
R8	& M01 F12 v M05 F10	& M02 F11 v M04 F05	& M03 F06 v M10 F02	&\\& M06 F09 v M12 F04	& M07 F08 v M08 F07	& M09 F01 v M11 F03 \\
R9	& M01 F03 v M08 F11	& M02 F07 v M11 F06	& M03 F02 v M12 F05	&\\& M04 F12 v M06 F08	& M05 F01 v M07 F04	& M09 F09 v M10 F10 \\
R10	& M01 F11 v M06 F04	& M02 F06 v M08 F03	& M03 F05 v M11 F07	&\\& M04 F02 v M07 F12	& M05 F08 v M12 F01	& M09 F10 v M10 F09 \\
R11	& M01 F10 v M04 F06	& M02 F08 v M05 F04	& M03 F01 v M08 F09	&\\& M06 F07 v M09 F02	& M07 F05 v M10 F03	& M11 F11 v M12 F12 \\
R12	& M01 F08 v M10 F05	& M02 F04 v M07 F10	& M03 F07 v M06 F01	&\\& M04 F09 v M08 F06	& M05 F02 v M09 F03	& M11 F12 v M12 F11 \\
\end{tabular}}
\end{Example}

\begin{Theorem}
A fully resolvable CMDRR$(3n, n)$ exists for all $n \ge 5$ and $n$ odd. 
\end{Theorem}
\begin{Proof}
Begin with an HSOLSSOM$(3^n)$ and convert to a resolved mixed doubles round robin tournament. Fill each hole with a 
CMDRR$(3,1)$ on the corresponding six players, noting that each CMDRR$(3,1)$ contributes one spouse pair to
the final schedule. Three of the four games from each CMDRR$(3,1)$ are added to the three rounds that
lack the six players from the hole. Collect together the fourth game from each CMDRR$(3,1)$ into one additional
round. This give a tournament with $3n+1$ rounds of play. The first $3n$ full rounds will all have 
$(3n-1)/2$ games and two byes, and the additional short round will have $n$ games and $2n$ byes. Over the course
of the tournament, each spouse pair player will receive exactly $2$ byes while the non-spouse pair players will 
receive exactly $1$ bye.\qed
\end{Proof}

\newpage
\begin{Example}\label{y} Lemma 2.2 of Abel et al.~\cite{ABZ} gives an example of a HSPS of type $3^5$
which is equivalent to the HSOLSSOM$(3^5)$ below.
\\[.05in]
{\footnotesize
\begin{tabular}{|ccccccccccccccc|}
\cline{1-15}
.  & .  & .  & 12 & 9  & 13 & 14 & 6  & 11 & 7  & 5  & 15 & 8  & 10 & 4  \\
.  & .  & .  & 14 & 10 & 7  & 12 & 15 & 4  & 13 & 8  & 6  & 5  & 9  & 11 \\
.  & .  & .  & 8  & 15 & 11 & 5  & 10 & 13 & 4  & 14 & 9  & 12 & 6  & 7  \\
7  & 11 & 13 & .  & .  & .  & 15 & 12 & 2  & 3  & 9  & 14 & 10 & 8  & 1  \\
14 & 8  & 12 & .  & .  & .  & 3  & 13 & 10 & 15 & 1  & 7  & 2  & 11 & 9  \\
10 & 15 & 9  & .  & .  & .  & 11 & 1  & 14 & 8  & 13 & 2  & 7  & 3  & 12 \\ 
4  & 13 & 11 & 10 & 14 & 2  & .  & .  & .  & 1  & 15 & 5  & 6  & 12 & 3  \\
12 & 5  & 14 & 3  & 11 & 15 & .  & .  & .  & 6  & 2  & 13 & 1  & 4  & 10 \\
15 & 10 & 6  & 13 & 1  & 12 & .  & .  & .  & 14 & 4  & 3  & 11 & 2  & 5  \\
6  & 9  & 15 & 7  & 2  & 14 & 13 & 3  & 5  & .  & .  & .  & 4  & 1  & 8  \\
13 & 4  & 7  & 15 & 8  & 3  & 6  & 14 & 1  & .  & .  & .  & 9  & 5  & 2  \\
8  & 14 & 5  & 1  & 13 & 9  & 2  & 4  & 15 & .  & .  & .  & 3  & 7  & 6  \\
11 & 7  & 4  & 9  & 12 & 1  & 10 & 5  & 3  & 2  & 6  & 8  & .  & .  & .  \\
5  & 12 & 8  & 2  & 7  & 10 & 1  & 11 & 6  & 9  & 3  & 4  & .  & .  & .  \\
9  & 6  & 10 & 11 & 3  & 8  & 4  & 2  & 12 & 5  & 7  & 1  & .  & .  & .  \\
\cline{1-15}
\end{tabular}
\\[.05in]
\begin{tabular}{|ccccccccccccccc|}
\cline{1-15}
.  & .  & .  & 14 & 10 & 7  & 12 & 15 & 4  & 13 & 8  & 6  & 5  & 9  & 11 \\
.  & .  & .  & 8  & 15 & 11 & 5  & 10 & 13 & 4  & 14 & 9  & 12 & 6  & 7  \\
.  & .  & .  & 12 & 9  & 13 & 14 & 6  & 11 & 7  & 5  & 15 & 8  & 10 & 4  \\
14 & 8  & 12 & .  & .  & .  & 3  & 13 & 10 & 15 & 1  & 7  & 2  & 11 & 9  \\
10 & 15 & 9  & .  & .  & .  & 11 & 1  & 14 & 8  & 13 & 2  & 7  & 3  & 12 \\ 
7  & 11 & 13 & .  & .  & .  & 15 & 12 & 2  & 3  & 9  & 14 & 10 & 8  & 1  \\
12 & 5  & 14 & 3  & 11 & 15 & .  & .  & .  & 6  & 2  & 13 & 1  & 4  & 10 \\
15 & 10 & 6  & 13 & 1  & 12 & .  & .  & .  & 14 & 4  & 3  & 11 & 2  & 5  \\
4  & 13 & 11 & 10 & 14 & 2  & .  & .  & .  & 1  & 15 & 5  & 6  & 12 & 3  \\
13 & 4  & 7  & 15 & 8  & 3  & 6  & 14 & 1  & .  & .  & .  & 9  & 5  & 2  \\
8  & 14 & 5  & 1  & 13 & 9  & 2  & 4  & 15 & .  & .  & .  & 3  & 7  & 6  \\
6  & 9  & 15 & 7  & 2  & 14 & 13 & 3  & 5  & .  & .  & .  & 4  & 1  & 8  \\
5  & 12 & 8  & 2  & 7  & 10 & 1  & 11 & 6  & 9  & 3  & 4  & .  & .  & .  \\
9  & 6  & 10 & 11 & 3  & 8  & 4  & 2  & 12 & 5  & 7  & 1  & .  & .  & .  \\
11 & 7  & 4  & 9  & 12 & 1  & 10 & 5  & 3  & 2  & 6  & 8  & .  & .  & .  \\
\cline{1-15}
\end{tabular}}
\\[.05in]
\newpage\noindent
The HSOLSSOM is converted to a mixed doubles tournament and filled to produce a CMDRR$(15,5)$ with $15$ full rounds and $1$ short round. 
The spouse pairs are M$1$F$1$, M$4$F$4$, M$7$F$7$, M$10$F$10$, and M$13$F$13$.
\\[.05in]
{\footnotesize
\begin{tabular}{ccccccccc}
R1	& M01 F02 v M02 F03	& M04 F09 v M11 F15	& M05 F13 v M08 F11	&\\& M06 F12 v M15 F08	& M07 F06 v M13 F10	& M09 F14 v M10 F05	&\\& M12 F07 v M14 F04 \\
R2	& M01 F03 v M03 F02	& M04 F10 v M13 F09	& M05 F07 v M12 F13	&\\& M06 F14 v M09 F12	& M07 F15 v M11 F06	& M08 F04 v M14 F11	&\\& M10 F08 v M15 F05 \\
R3	& M02 F01 v M03 F03	& M04 F15 v M07 F10	& M05 F11 v M14 F07	&\\& M06 F08 v M10 F14	& M08 F13 v M12 F04	& M09 F05 v M15 F12	&\\& M11 F09 v M13 F06 \\
R4	& M01 F11 v M09 F15	& M02 F13 v M10 F09	& M03 F07 v M15 F10	&\\& M04 F05 v M05 F06	& M07 F12 v M14 F01	& M08 F02 v M11 F14	&\\& M12 F03 v M13 F08 \\
R5	& M01 F08 v M13 F11	& M02 F12 v M07 F13	& M03 F14 v M11 F07	&\\& M04 F06 v M06 F05	& M08 F10 v M15 F02	& M09 F03 v M12 F15	&\\& M10 F01 v M14 F09 \\
R6	& M01 F15 v M12 F08	& M02 F09 v M14 F12	& M03 F10 v M08 F14	&\\& M05 F04 v M06 F06	& M07 F01 v M10 F13	& M09 F11 v M13 F03	&\\& M11 F02 v M15 F07 \\
R7	& M01 F13 v M06 F10	& M02 F11 v M15 F06	& M03 F04 v M10 F15	&\\& M04 F14 v M12 F01	& M05 F02 v M13 F12	& M07 F08 v M08 F09	&\\& M11 F05 v M14 F03 \\
R8	& M01 F05 v M11 F13	& M02 F14 v M04 F11	& M03 F12 v M13 F04	&\\& M05 F15 v M10 F02	& M06 F03 v M14 F10	& M07 F09 v M09 F08	&\\& M12 F06 v M15 F01 \\
R9	& M01 F10 v M14 F05	& M02 F06 v M12 F14	& M03 F15 v M05 F12	&\\& M04 F01 v M15 F11	& M06 F13 v M11 F03	& M08 F07 v M09 F09	&\\& M10 F04 v M13 F02 \\
R10	& M01 F09 v M05 F14	& M02 F15 v M08 F05	& M03 F06 v M14 F08	&\\& M04 F02 v M09 F13	& M06 F07 v M13 F01	& M07 F03 v M15 F04	&\\& M10 F11 v M11 F12 \\
R11	& M01 F04 v M15 F09	& M02 F07 v M06 F15	& M03 F13 v M09 F06	&\\& M04 F08 v M14 F02	& M05 F03 v M07 F14	& M08 F01 v M13 F05	&\\& M10 F12 v M12 F11 \\
R12	& M01 F14 v M07 F04	& M02 F05 v M13 F07	& M03 F08 v M04 F13	&\\& M05 F09 v M15 F03	& M06 F01 v M08 F15	& M09 F02 v M14 F06	&\\& M11 F10 v M12 F12 \\
R13	& M01 F07 v M10 F06	& M02 F04 v M09 F10	& M03 F11 v M06 F09	&\\& M04 F12 v M08 F03	& M05 F01 v M11 F08	& M07 F05 v M12 F02	&\\& M13 F14 v M14 F15 \\
R14	& M01 F12 v M04 F07	& M02 F08 v M11 F04	& M03 F05 v M07 F11	&\\& M05 F10 v M09 F01	& M06 F02 v M12 F09	& M08 F06 v M10 F03	&\\& M13 F15 v M15 F14 \\
R15	& M01 F06 v M08 F12	& M02 F10 v M05 F08	& M03 F09 v M12 F05	&\\& M04 F03 v M10 F07	& M06 F11 v M07 F02	& M09 F04 v M11 F01	&\\& M14 F13 v M15 F15 \\
R16	& M02 F02 v M03 F01	& M05 F05 v M06 F04	& M08 F08 v M09 F07	&\\& M11 F11 v M12 F10	& M14 F14 v M15 F13 \\
\end{tabular}}
\end{Example}

\newpage
\begin{Example}\label{C166} Abel et al.~\cite{ABZZ} gives an example of an HSOLS$(3^51^1)$ \\ which is shown in 
block diagonal form below.  
\\[.05in]
{\scriptsize
\begin{tabular}{|cccccccccccccccc|}
\cline{1-16}
 .$\!$ & .$\!$  & .$\!$ & 10$\!$ & 14$\!$ & 11$\!$ & 5$\!$  & 12$\!$ & 13$\!$ & 15$\!$ & 7$\!$  & 8$\!$  & 6$\!$  & 16$\!$ & 4$\!$  & 9  \\
 .$\!$  & .$\!$  & .$\!$  & 12$\!$ & 13$\!$ & 14$\!$ & 10$\!$ & 16$\!$ & 4$\!$  & 9$\!$  & 15$\!$ & 5$\!$  & 11$\!$ & 6$\!$  & 8$\!$  & 7  \\
 .$\!$  & .$\!$  & .$\!$  & 13$\!$ & 15$\!$ & 10$\!$ & 4$\!$  & 11$\!$ & 12$\!$ & 14$\!$ & 6 $\!$ & 9 $\!$ & 7 $\!$ & 5 $\!$ & 16$\!$ & 8  \\
 9 $\!$ & 16$\!$ & 11$\!$ & . $\!$ & . $\!$ & . $\!$ & 14$\!$ & 3 $\!$ & 2 $\!$ & 13$\!$ & 8 $\!$ & 15$\!$ & 10$\!$ & 1 $\!$ & 7 $\!$ & 12 \\
16 $\!$ & 8 $\!$ & 7 $\!$ & . $\!$ & . $\!$ & . $\!$ & 15$\!$ & 14$\!$ & 1 $\!$ & 3 $\!$ & 9 $\!$ & 13$\!$ & 12$\!$ & 2 $\!$ & 11$\!$ & 10 \\
 7 $\!$ & 12$\!$ & 15$\!$ & . $\!$ & . $\!$ & . $\!$ & 1 $\!$ & 13$\!$ & 3 $\!$ & 8 $\!$ & 14$\!$ & 16$\!$ & 9 $\!$ & 10$\!$ & 2 $\!$ & 11 \\
14 $\!$ & 6 $\!$ & 10$\!$ & 2 $\!$ &11 $\!$ & 13$\!$ & . $\!$ & . $\!$ & . $\!$ & 4 $\!$ & 5 $\!$ & 3 $\!$ & 16$\!$ & 12$\!$ & 1 $\!$ & 15 \\
 5 $\!$ & 11$\!$ & 14$\!$ & 15$\!$ & 10$\!$ & 16$\!$ & . $\!$ & . $\!$ & . $\!$ & 6 $\!$ & 1 $\!$ & 4 $\!$ & 2 $\!$ & 3 $\!$ & 12$\!$ & 13 \\
 6 $\!$ & 15$\!$ & 13$\!$ & 11$\!$ & 16$\!$ & 12$\!$ & . $\!$ & . $\!$ & . $\!$ & 5 $\!$ & 2 $\!$ & 1 $\!$ & 3 $\!$ & 4 $\!$ & 10$\!$ & 14 \\
 8 $\!$ & 14$\!$ & 4 $\!$ & 3 $\!$ & 7 $\!$ & 2 $\!$ & 13$\!$ & 15$\!$ & 16$\!$ & . $\!$ & . $\!$ & . $\!$ & 5 $\!$ & 9 $\!$ & 6 $\!$ & 1  \\
13 $\!$ & 5 $\!$ & 16$\!$ & 14$\!$ & 1 $\!$ & 7 $\!$ & 3 $\!$ & 6 $\!$ & 15$\!$ & . $\!$ & . $\!$ & . $\!$ & 4 $\!$ & 8 $\!$ & 9 $\!$ & 2  \\
 4 $\!$ & 13$\!$ & 6 $\!$ & 1 $\!$ & 9 $\!$ & 15$\!$ & 16$\!$ & 2 $\!$ & 14$\!$ & . $\!$ & . $\!$ & . $\!$ & 8 $\!$ & 7 $\!$ & 5 $\!$ & 3  \\
11 $\!$ & 9 $\!$ & 12$\!$ & 7 $\!$ & 8 $\!$ & 3 $\!$ & 2 $\!$ & 5 $\!$ & 10$\!$ & 1 $\!$ & 16$\!$ & 6 $\!$ & . $\!$ & . $\!$ & . $\!$ & 4  \\
10 $\!$ & 7 $\!$ & 9 $\!$ & 8 $\!$ & 12$\!$ & 1 $\!$ & 6 $\!$ & 4 $\!$ & 11$\!$ & 16$\!$ & 3 $\!$ & 2 $\!$ & . $\!$ & . $\!$ & . $\!$ & 5  \\
12 $\!$ & 10$\!$ & 8 $\!$ & 16$\!$ & 3 $\!$ & 9 $\!$ & 11$\!$ & 1 $\!$ & 5 $\!$ & 2 $\!$ & 4 $\!$ & 7 $\!$ & . $\!$ & . $\!$ & . $\!$ & 6  \\
15 $\!$ & 4 $\!$ & 5 $\!$ & 9 $\!$ & 2 $\!$ & 8 $\!$ & 12$\!$ & 10$\!$ & 6 $\!$ & 7 $\!$ & 13$\!$ & 14$\!$ & 1 $\!$ & 11$\!$ & 3 $\!$ & .  \\
\cline{1-16}
\end{tabular}}
\\[.05in]
A CMDRR$(16,6)$ can be derived from this and can be played in 25 short rounds of 5 games each (by computer search).
The spouse pairs are M$1$F$1$,  M$4$F$4$,  M$7$F$7$,  M$10$F$10$,  M$13$F$13$, and  M$16$F$16$.
\\[.05in]
{\footnotesize
\begin{tabular}{ccccccc}
R1  & M02 F11 v M13 F09	& M04 F08 v M11 F14	& M05 F05 v M06 F04	&\\& M07 F15 v M16 F12	& M10 F06 v M15 F02 \\
R2	& M03 F04 v M07 F10	& M05 F12 v M13 F08	& M06 F16 v M12 F15	&\\& M08 F01 v M11 F06	& M14 F14 v M15 F13 \\
R3	& M01 F07 v M11 F13	& M02 F08 v M15 F10	& M03 F14 v M10 F04	&\\& M04 F02 v M09 F11	& M07 F12 v M14 F06 \\
R4	& M01 F09 v M16 F15	& M04 F13 v M10 F03	& M05 F02 v M14 F12	&\\& M06 F14 v M11 F07	& M09 F10 v M15 F05 \\
R5	& M02 F04 v M09 F15	& M03 F11 v M08 F14	& M06 F01 v M07 F13	&\\& M10 F09 v M14 F16	& M11 F10 v M12 F12 \\
R6	& M02 F10 v M07 F06	& M04 F03 v M08 F15	& M05 F01 v M09 F16	&\\& M10 F11 v M11 F12	& M12 F05 v M15 F07 \\
R7	& M01 F13 v M09 F06	& M02 F09 v M10 F14	& M03 F16 v M15 F08	&\\& M07 F05 v M11 F03	& M08 F04 v M12 F02 \\
R8	& M01 F02 v M02 F03	& M03 F07 v M13 F12	& M04 F06 v M06 F05	&\\& M08 F13 v M16 F10	& M09 F01 v M12 F14 \\
R9	& M04 F15 v M12 F01	& M07 F04 v M10 F13	& M08 F02 v M13 F05	&\\& M09 F14 v M16 F06	& M11 F08 v M14 F03 \\
R10	& M05 F14 v M08 F10	& M06 F02 v M15 F09	& M07 F03 v M12 F16	&\\& M09 F04 v M14 F11	& M10 F05 v M13 F01 \\
R11	& M02 F05 v M12 F13	& M03 F15 v M05 F07	& M06 F10 v M14 F01	&\\& M11 F04 v M13 F16	& M15 F06 v M16 F03 \\
R12	& M02 F14 v M06 F12	& M03 F09 v M12 F06	& M04 F01 v M14 F08	&\\& M05 F03 v M10 F07	& M11 F02 v M16 F13 \\
R13	& M01 F08 v M12 F04	& M05 F09 v M11 F01	& M06 F03 v M09 F12	&\\& M07 F16 v M13 F02	& M08 F06 v M10 F15 \\
R14	& M01 F11 v M06 F07	& M02 F02 v M03 F01	& M08 F03 v M14 F04	&\\& M09 F05 v M10 F16	& M13 F15 v M15 F14 \\

\end{tabular}
\\[.05in]
\begin{tabular}{ccccccc}

R15	& M01 F04 v M15 F12	& M03 F08 v M16 F05	& M08 F07 v M09 F09	&\\& M11 F11 v M12 F10	& M13 F14 v M14 F15 \\
R16	& M01 F10 v M04 F09	& M03 F12 v M09 F13	& M05 F15 v M07 F11	&\\& M10 F01 v M16 F07	& M12 F08 v M13 F06 \\
R17	& M01 F15 v M10 F08	& M02 F06 v M14 F07	& M05 F10 v M16 F02	&\\& M06 F13 v M08 F16	& M11 F09 v M15 F04 \\ 
R18	& M01 F14 v M05 F16	& M04 F10 v M13 F07	& M07 F08 v M08 F09	&\\& M09 F02 v M11 F15	& M14 F05 v M16 F11 \\
R19	& M01 F16 v M14 F10	& M02 F15 v M11 F05	& M03 F13 v M04 F11	&\\& M05 F04 v M06 F06	& M12 F03 v M16 F14 \\
R20	& M01 F05 v M07 F14	& M02 F13 v M05 F08	& M04 F07 v M15 F16	&\\& M06 F09 v M13 F03	& M10 F12 v M12 F11 \\
R21	& M03 F06 v M11 F16	& M04 F12 v M16 F09	& M07 F01 v M15 F11	&\\& M09 F03 v M13 F10	& M12 F07 v M14 F02 \\
R22	& M01 F03 v M03 F02	& M02 F16 v M08 F11	& M04 F05 v M05 F06	&\\& M07 F09 v M09 F08	& M14 F13 v M15 F15 \\
R23	& M01 F06 v M13 F11	& M02 F07 v M16 F04	& M03 F05 v M14 F09	&\\& M06 F08 v M10 F02	& M08 F12 v M15 F01 \\
R24	& M01 F12 v M08 F05	& M02 F01 v M03 F03	& M04 F14 v M07 F02	&\\& M05 F13 v M12 F09	& M06 F11 v M16 F08 \\
R25	& M02 F12 v M04 F16	& M03 F10 v M06 F15	& M05 F11 v M15 F03	&\\& M08 F08 v M09 F07	& M13 F04 v M16 F01 \\
\end{tabular}}
\end{Example}

   \section{Product Theorem}

We next present a product construction for CMDRR. While this does not expand the spectrum given in 
Section 3, it does provide an alternative construction that does not rely on HSOLS.

\begin{Theorem}\label{product}
If there exists a CMDRR$(n,k)$, a SAMDRR$(m)$, and two mutually orthogonal latin squares, MOLS, of order $n$, 
then there exists a  CMDRR$(mn,mk)$.
\end{Theorem}
\begin{Proof}
Let M$(i)$ and F$(i)$, with $i=1, \ldots, n$, denote the players of the CMDRR$(n,k)$, and
let M$'(j)$ and F$'(j)$, with $j=1, \ldots, m$, denote the players of the SAMDRR$(m)$. As usual
assume, without loss of generality, that spouses have the same index.
We will construct a CMDRR$(mn,mk)$ on new players M$(i,j)$ and F$(i,j)$, with $i=1, \ldots, n$,
and $j=1, \ldots, m$. 

For each game M$(w)$F$(x)$ v M$(y)$F$(z)$ of the CMDRR$(n,k)$, add to the new CMDRR$(mn,mk)$ the $m$ games 
 M$(w,j)$F$(x,j)$ v M$(y,j)$F$(z,j)$, with $j=1, \ldots, m$. Call these type $1$ games. There are $m(n^2-k)/2$ 
of these games. 

Let L$_1$ and L$_2$ be the two MOLS of order $n$. 
For each game M$'(w)$F$'(x)$ v M$'(y)$F$'(z)$ of the SAMDRR$(m)$, add to the new CMDRR$(mn,mk)$ the $n^2$ games 
 M$(i_1,w)$F$(i_2,x)$ v M$(i_3,y)$F$(i_4,z)$, with $i_1,i_2=1, \ldots, n$, and
$i_3=$L$_1(i_1,i_2)$, and $i_4=$L$_2(i_1,i_2)$.  Note that all of $w,x,y,$ and $z$ are distinct. Call these type $2$ games.
There are $n^2(m^2-m)/2$ of these games. So the total number of type 1 and type 2 games 
is $m(n^2-k)/2 + n^2(m^2-m)/2 = ((mn)^2 - mk)/2$, the number of games expected for a CMDRR$(mn,mk)$.

We now check that the conditions for a CMDRR$(mn,mk)$ are met for opposite sex players. If M$(i)$ and F$(i)$ are spouses in the 
CMDRR$(n,k)$, then for each $j=1, \ldots, m$, the players M$(i,j)$ and F$(i,j)$ satisfy the condition for spouses in the CMDRR$(mn,mk)$, 
because each pair never occurs in a type $1$ or type $2$ game as partners or opponents. Thus there are at least $mk$
spouse pairs. Consider any other pair  M$(i_1,w)$F$(i_2,x)$ that are not one of these spouse pairs. If $w=x$ then by 
construction the players partner once and oppose once in type 1 games. If $w\ne x$ then M$'(w)$ and F$'(x)$ partner
and oppose exactly once in the SAMDRR$(m)$ and by definition of MOLS, M$(i_1,w)$ and F$(i_2,x)$ 
partner and oppose exactly once in the CMDRR$(mn,mk)$. We conclude that there are exactly $mk$ spouse pairs and that every male
and female who are not spouses are partners exactly once and opponents exactly once.

We now check that the conditions for a CMDRR$(mn,mk)$ are met for same sex players. Consider players 
M$(i_1,w)$ and M$(i_3,y)$. If $w=y$ then by construction they oppose at least once
in a type 1 game. If $w\ne y$ then again by construction they oppose exactly once in a type 2 game. The condition for 
female players is analogous. So same sex players oppose at least once. The total number of games is correct so we conclude
that each player who does not have a spouse opposes 
some other same sex player who does not have a spouse exactly twice and opposes all other same sex players exactly once. \qed
\end{Proof}



\begin{thebibliography}{77}
\frenchspacing

\bibitem{ABZ} R.J.R.\ Abel, F.E.\ Bennett, and H.\ Zhang, Holey Steiner pentagon systems, J. Combin. Des. {\bf 7} (1999), 41--56. 

\bibitem{ABZZ} R.J.R.\ Abel, F.E.\ Bennett, H.\ Zhang, and L.\ Zhu, A few more incomplete self-orthogonal latin squares
and related designs, Australas. J. Combin. {\bf 21} (2000), 85--94. 

\bibitem{IA1} I.\ Anderson, {\it Combinatorial Designs and Tournaments}, Oxford University Press, Oxford, 1997.

\bibitem{IA2} I.\ Anderson, Early examples of spouse avoidance, Bull. Inst. Combin. Appl. {\bf 54} (2008),
47--52.

\bibitem{BZ} F.E.\ Bennett and L.\ Zhu, The spectrum of HSOLSSOM$(h^n)$ where $h$ is even, Discrete Math. {\bf 158} (1996), 11--25. 

\bibitem{BS1}	D.R.\ Berman and D.D.\ Smith, Mitchell tournaments, Bull. Inst. Combin. Appl. {\bf 65} (2012), 33--42.

\bibitem{BS2}	D.R.\ Berman and D.D.\ Smith, Balanced equitable mixed doubles round-robin tournaments, 
Bull. Inst. Combin. Appl. {\bf 68} (2013), 90--101.

\bibitem{BCH1} R.K.\ Brayton, D.\ Coppersmith, and A.J.\ Hoffman,
Self-orthogonal latin squares of all orders $n \ne 2, 3, 6$, Bull.
Amer. Math. Soc. {\bf 80} (1974) 116--118.

\bibitem{BCH2} R.K.\ Brayton, D.\ Coppersmith, and A.J.\ Hoffman,
Self-orthogonal latin squares, Coll. Int. Th. Comb., Rome 1973, Atti
del Convegni Lincei, {\bf 17} (1976) 509--517.

\bibitem{FIN} N.J.\ Finizio and L.\ Zhu,
Self-orthogonal latin squares, in: C.J. Colbourn and J.H. Dinitz
(eds.), {\em Handbook of Combinatorial Designs}, second edition, CRC
Press, Boca Raton, FL, 2007, 211--219.

\bibitem{GA} G.\ Ge and R.J.R.\ Abel, Some new HSOLSSOMs of type $h^n$ and $1^nu^1$, J. Combin. Des.
{\bf 9} (2001), 435--444. 

\bibitem{HZ} K.\ Heinrich and L.\ Zhu, Incomplete self-orthogonal latin squares, 
J. Austral. Math. Soc. (Series A) {\bf 42} (1987), 365--384.

\bibitem{M} J.T.\ Mitchell, {\it Duplicate Whist}, 2nd edition,
Ihling Bros., Kalamazoo, 1897.

\bibitem{Mor} L.B.\ Morales, Constructing cyclic PBIBD$(2)$s through an optimization approach: thirty-two new cyclic designs, 
J. Combin. Des. {\bf 13} (2005), 377--387. 

\bibitem{XZZ} Y.\ Xu, H.\ Zhang, and L.\ Zhu, Existence of frame SOLS of type $a^nb^1$, Discrete
Math. {\bf 250} (2002), 211--230.

\bibitem{Z} L.\ Zhu, Existence of self-orthogonal latin squares ISOLS$(6m+2,2m)$, Ars Combin. 
{\bf 39} (1995), 65--74.

\nonfrenchspacing
\end{thebibliography}
\end{document}